\documentclass[10pt,sumlimits,namelimits]{article}

\usepackage{amsthm}
\usepackage{amsmath}
\usepackage{amsfonts}
\usepackage{amssymb}
\usepackage{mathrsfs}
\usepackage{stmaryrd}
\usepackage{dsfont}
\usepackage{enumerate}
\usepackage{color}
\usepackage{hyperref}
\usepackage[applemac]{inputenc}

\newcommand{\A}{\mathbb{A}}

\newcommand{\N}{\mathbb{N}}
\newcommand{\Z}{\mathbb{Z}}

\newcommand{\R}{\mathbb{R}}
\newcommand{\C}{\mathbb{C}}
\renewcommand{\O}{\Omega}

\renewcommand{\S}{\mathbb{S}}

\newcommand{\vO}{{\cal O}}

\newcommand{\vR}{{\cal R}}

\newcommand{\Dr}{\mathscr{D}}

\newcommand{\vphi}{\varphi}
\newcommand{\eps}{\varepsilon}
\newcommand{\la}{\lambda}
\newcommand{\dsp}{\displaystyle}
\newcommand{\ovl}{\overline}
\newcommand{\udl}{\underline}
\newcommand{\vlim}{\lim\limits}

\newcommand{\vmin}{\min\limits}

\newcommand{\vint}{\int\limits}
\newcommand{\vsum}{\sum\limits}
\newcommand{\inj}{\hookrightarrow}
\newcommand{\tends}{\longrightarrow}

\newcommand{\wt}{\widetilde}

\newcommand{\li}{\llbracket}
\newcommand{\ri}{\rrbracket}
\newcommand{\loc}{\mathrm{loc}}

\renewcommand{\c}{\mathrm{c}}
\renewcommand{\d}{\mathrm{d}}

\newcommand{\dist}{\mathrm{dist}}

\newcommand{\M}{\mathrm{max}}

\renewcommand{\le}{\leqslant}
\renewcommand{\ge}{\geqslant}
\renewcommand{\Re}{\mathrm{Re}}
\renewcommand{\Im}{\mathrm{Im}}

\newcommand{\vi}{\mathrm{i}}
\newcommand{\p}{\prime}

\DeclareMathOperator{\supp}{supp}

\numberwithin{equation}{section}

\newtheorem{thm}{Theorem}[section]

\newtheorem{cor}[thm]{Corollary}
\newtheorem{lem}[thm]{Lemma}
\newtheorem{sym}[thm]{Symmetry Property}

\theoremstyle{definition}
\newtheorem{rmk}[thm]{Remark}

\newtheorem{exa}[thm]{Example}

\newtheorem{ass}[thm]{Assumption}

\newenvironment{proof*}{\noindent{\bf Proof.}}{\qed}
\newenvironment{vproof}[1]{\noindent{\bf Proof #1}}{\qed}

\title{\huge \sc A sharper energy method for the localization of the support to some stationary Schrödinger equations with a singular nonlinearity}
\author{\sc Pascal Bégout$^*$ and Jes\'us Ildefonso D\'iaz$^\dagger$}
\date{}

\begin{document}

\maketitle

\begin{gather*}
\begin{array}{cc}
            ^*\mbox{Institut de Mathématiques de Toulouse \& TSE }	&	\;^\dagger\mbox{Instituto de Matem\'atica Interdisciplinar}		\\
                                          \mbox{Université Toulouse I Capitole }	&	\mbox{ Departamento de Matem\'atica Aplicada}			\\
                                                   \mbox{Manufacture des Tabacs }	&	\mbox{ Universidad Complutense de Madrid}				\\
                                                          \mbox{21, Allée de Brienne }	&	\mbox{ Plaza de las Ciencias, 3}						\\
                                  \mbox{31015 Toulouse Cedex 6, FRANCE }	&	\mbox{ 28040 Madrid, SPAIN}
\bigskip \\
\mbox{
{\footnotesize $^*$e-mail\:: \href{mailto:Pascal.Begout@math.cnrs.fr}{\udl{\texttt{Pascal.Begout@math.cnrs.fr}}}}}
&
\mbox{
{\footnotesize $^\dagger$e-mail\:: \href{mailto:diaz.racefyn@insde.es}{\udl{\texttt{diaz.racefyn@insde.es}}}}
}
\end{array}
\end{gather*}

\begin{abstract}
We prove the compactness of the support of the solution of some stationary Schrödinger equations with a singular nonlinear order term. We present here a sharper version of some energy methods previously used in the literature and, in particular, by the authors.
\end{abstract}

{\let\thefootnote\relax\footnotetext{$^\dagger$The research of J.I.~D\'iaz was partially supported by the project ref. MTM200806208 of the DGISPI (Spain) and the Research Group MOMAT (Ref. 910480) supported by UCM. He has received also support from the ITN \textit{FIRST} of the Seventh Framework Program of the European Community's (grant agreement number 238702)}}
{\let\thefootnote\relax\footnotetext{2010 Mathematics Subject Classification: 35B45}}
{\let\thefootnote\relax\footnotetext{Key Words: energy method, Schrödinger equation, solutions with compact support}}

\tableofcontents

\baselineskip .65cm

\section{Introduction}
\label{intro}

Since the beginnings of the eighties of the last century, it is already well-known that the absence of the maximum principle for the case of systems and higher order nonlinear partial differential equations was one of the main motivations of the introduction of suitable energy methods allowing to conclude the compactness of the support of their solutions (see, e.g., the presentation made in the monograph Antontsev, D\'iaz and
Shmarev~\cite{MR2002i:35001}).

\medskip
\noindent
The application of such type of methods to the case of nonlinear Schrödinger equations with a singular zero order term required some important improvements of the method. That was the main object of the previous author's papers of Bégout and D\'iaz~\cite{MR2214595,MR2876246}.

\medskip
\noindent
The main goal of this new paper is to present a sharper version of the mentioned method potentially able to be applied to many other problems related to this type of Schrödinger equations such as the study of self-similar solutions, case of Neumann boundary conditions, presence of nonlocal terms (such as, for instance, in Hartree-Fock theory: Cazenave~\cite{MR2002047}), etc., which can not be treated with the mere technique presented in Bégout and D\'iaz~\cite{MR2214595,MR2876246}. As a matter of fact, the concrete application of this sharper energy method to the concrete case of self-similar solutions of the evolution Schrödinger problem requires many additional arguments justifying the special structure of those solutions, reason why we decided to present it in a separated work (Bégout and D\'iaz~\cite{MR3193996}). We send the reader to Bégout and D\'iaz~\cite{MR3193996} for a long description of the important role of the compactness of the solution in this context and for many other references related to this qualitative property of the solution.

\medskip
\noindent
This paper is organized as follows. Below, we give some notations which will be used throughout this paper. In Section~\ref{grl}, we give the precise ``localization'' estimates which imply a solution of a partial differential equation to be compactly supported \Big(see Theorems~\ref{thmsta} and \ref{thmstaF}, and especially estimates~\eqref{thmsta1} and \eqref{thmstaF1}\Big). In Section~\ref{gfs}, we give a tool which permits, from a solution of some partial differential equation, to establish the ``localization'' estimate (Theorem~\ref{thmestgen}). The results of these two sections are proved in Section~\ref{proof}. In Bégout and D\'iaz~\cite{MR2876246}, localization property is studied for the complex-valued equation
\begin{gather}
\label{nls}
-\Delta u+a|u|^{-(1-m)}u+bu=F, \text{ in } \O.
\end{gather}
We also study this property here, but with a change of notation (see Remark~\ref{com} below for the motivation of this change). Section~\ref{app} is devoted to the study of the localization property of the solutions of equation~\eqref{nls}, in the same spirit as Bégout and D\'iaz~\cite{MR2876246}, but with the homogeneous Neumann boundary condition instead of the homogeneous Dirichlet boundary condition (compare Theorem~\ref{thmcomneu1} below with Theorem~3.5 in Bégout and D\'iaz~\cite{MR2876246}). Finally, at the end of the paper, we treat equation~\eqref{nls} with the homogeneous Dirichlet boundary condition (Remark~\ref{rmkfinal}). We state the same results as in Bégout and D\'iaz~\cite{MR2876246}, but with now the weaker assumption $F\in L^2(\O).$

\medskip
\noindent
Before ending this section, we shall indicate here some of the notations used throughout. We write $\vi^2=-1.$ We denote by $\ovl z$ the conjugate of the complex number $z.$ For $1\le p\le\infty,$ $p^\p$ is the conjugate of $p$ defined by
$\frac1p+\frac1{p^\p}=1.$ For $j,k\in\Z$ with $j<k,$ $\li j,k\ri=[j,k]\cap\Z.$ We denote by $\Gamma$ the boundary of a nonempty subset $\O\subseteq\R^N$ and $\O^\c=\R^N\setminus\O$ its complement. Unless if specified, any function lying in a functional space $\big(L^p(\O),$ $W^{m,p}(\O),$ etc\big) is supposed to be a complex-valued function $\big(L^p(\O;\C),$ $W^{m,p}(\O;\C),$ etc\big). For a Banach space $E,$ we denote by $E^\star$ its topological dual and by $\langle\: . \; , \: . \:\rangle_{E^\star,E}\in\R$ the $E^\star-E$ duality product. In particular, for any $T\in L^{p^\p}(\O)$ and $\vphi\in L^p(\O)$ with $1\le p<\infty,$
$\langle T,\vphi\rangle_{L^{p^\p}(\O),L^p(\O)}=\Re\vint_\O T(x)\ovl{\vphi(x)}\d x.$ As usual, we denote by $C$ auxiliary positive constants, and sometimes, for positive parameters $a_1,\ldots,a_n,$ write $C(a_1,\ldots,a_n)$ to indicate that the constant $C$ continuously depends only on $a_1,\ldots,a_n$ (this convention also holds for constants which are not denoted by ``$C$'').

\section{From suitable local inequalities to the vanishing of the involved complex functions on some small ball}
\label{grl}

In this section, we establish some results improving the presentation of some energy methods of Antontsev, D\'iaz and Shmarev~\cite{MR2002i:35001} which allow to prove localization properties of solutions of a general class of nonlinear partial differential equations (Section~\ref{app},
Remark~\ref{rmkfinal} below and Bégout and D\'iaz~\cite{MR3193996}). In contrast to the presentation in Bégout and D\'iaz~\cite{MR2876246} (see e.g. Theorem~1.1), the following statement does not need any information on the second order equation but it will merely use a suitable balance between the total local energy (diffusion + absorption local energies) and the local boundary flux. This will be crucial for the applicability of the method to cases for which the techniques of Bégout and D\'iaz~\cite{MR2214595,MR2876246} can not be applied.

\begin{thm}
\label{thmsta}
Assume $0<m<1$ and let $N\in\N.$ Then there exists $C=C(N,m)$ satisfying the following property$:$ let $x_0\in\R^N,$ $\rho_0>0$ and
$u\in H^1_\loc\big(B(x_0,\rho_0)\big).$ If there exist $L>0$ and $M>0$ such that for almost every $\rho\in(0,\rho_0),$
\begin{gather}
\label{thmsta1}
\|\nabla u\|_{L^2(B(x_0,\rho))}^2+L\|u\|_{L^{m+1}(B(x_0,\rho))}^{m+1}\le M\left|\dsp\int_{\S(x_0,\rho)}u\ovl{\nabla u}.\frac{x-x_0}{|x-x_0|}\d\sigma\right|,
\end{gather}
then $u_{|B(x_0,\rho_\M)}\equiv0,$ where
\begin{multline}
 \label{thmsta2}
  \rho_\M^\nu=\left(\rho_0^\nu-CM^2\max\left\{1,\frac{1}{L^2}\right\}\max\left\{\rho_0^{\nu-1},1\right\}\right. \\
  \left.\times\vmin_{\tau\in\left(\frac{m+1}{2},1\right]}\left\{\frac{E(\rho_0)^{\gamma(\tau)}
  \max\{b(\rho_0)^{\mu(\tau)},b(\rho_0)^{\eta(\tau)}\}}{2\tau-(1+m)}\right\}\right)_+,
\end{multline}
and where,
\begin{gather*}
\begin{array}{lll}
E(\rho_0)=\|\nabla u \|_{L^2(B(x_0,\rho_0))}^2,		&	&	b(\rho_0)=\|u\|_{L^{m+1}(B(x_0,\rho_0))}^{m+1},	\medskip \\
k=2(1+m)+N(1-m),							&	&	\nu=\frac{k}{m+1}>2,
\end{array}
\\
\gamma(\tau)=\frac{2\tau-(1+m)}{k}\in(0,1), \quad \mu(\tau)=\frac{2(1-\tau)}{k}, \quad \eta(\tau)=\frac{1-m}{1+m}-\gamma(\tau)>0.
\end{gather*}
for any $\tau\in\left(\frac{m+1}{2},1\right].$
\end{thm}

\noindent
Here and in what follows, $r_+=\max\{0,r\}$ denotes the positive part of the real number $r.$ For $x_0\in\R^N$ and $r>0,$ $B(x_0,r)$ is the open ball of $\R^N$ of center $x_0$ and radius $r,$ $\S(x_0,r)$ is its boundary and $\ovl B(x_0,r)$ is its closure. Finally, $\sigma$ is the surface measure on a sphere. A sharper estimate, in the same line of extension of the applicability of the techniques of Bégout and D\'iaz~\cite{MR2214595,MR2876246} indicated before, can be obtained under some additional assumption on $F.$

\begin{thm}
\label{thmstaF}
Let $0<m<1,$ $x_0\in\R^N,$ $\rho_1>\rho_0>0,$ $F\in L^2\big(B(x_0,\rho_1)\big)$ and $u\in H^1_\loc\big(B(x_0,\rho_1)\big).$ If there exist $L>0$ and $M>0$ such that for almost every $\rho\in(0,\rho_1),$
\begin{multline}
\label{thmstaF1}
\|\nabla u\|_{L^2(B(x_0,\rho))}^2+L\|u\|_{L^{m+1}(B(x_0,\rho))}^{m+1}+L\|u\|_{L^2(B(x_0,\rho))}^2		\\
\le M\left(\left|\int_{\S(x_0,\rho)}u\ovl{\nabla u}.\frac{x-x_0}{|x-x_0|}\d\sigma\right|+\int_{B(x_0,\rho)}|F(x)u(x)|\d x\right),
\end{multline}
then there exist $E_\star>0$ and $\eps_\star>0$ satisfying the following property$:$ if $\|\nabla u\|_{L^2(B(x_0,\rho_1))}^2<E_\star$ and
\begin{gather}
\label{thmstaF2}
\|F\|_{L^2(B(x_0,\rho))}^2\le\eps_\star\big((\rho-\rho_0)_+\big)^p, \; \forall\rho\in(0,\rho_1),
\end{gather}
where $p=\frac{2(1+m)+N(1-m)}{1-m},$ then $u_{|B(x_0,\rho_0)}\equiv0.$ In other words, with the notation of Theorem~$\ref{thmsta},$ $\rho_\M=\rho_0.$
\end{thm}

\begin{rmk}
\label{rmkthmstaF}
We may estimate $E_\star$ and $\eps_\star$ as
\begin{gather*}
E_\star=E_\star\left(\|u\|_{L^{m+1}(B(x_0,\rho_1))}^{-1},\rho_1,\frac{\rho_0}{\rho_1},\frac{L}{M},N,m\right), \\
\eps_\star=\eps_\star\left(\|u\|_{L^{m+1}(B(x_0,\rho_1))}^{-1},\frac{\rho_0}{\rho_1},\frac{L}{M},N,m\right).
\end{gather*}
The dependence on $\frac{1}{\delta}$ means that if $\delta$ goes to $0$ then $E_\star$ and $\eps_\star$ may be very large. Note that
$p=\frac{1}{\gamma(1)},$ where $\gamma$ is the function defined in Theorem~\ref{thmsta}.
\end{rmk}

\begin{rmk}
\label{rmkthmsta}
Note that by Cauchy-Schwarz's inequality, the right-hand side in \eqref{thmsta1} belongs to $L^1_\loc([0,\rho_0);\R)$ and so is defined almost everywhere in $(0,\rho_0).$ Consequently, by Hölder's inequality, the right-hand side in \eqref{thmstaF1} is defined almost everywhere in $(0,\rho_1).$
\end{rmk}

\section{A general framework of applications related to the Schrödinger operator}
\label{gfs}

The following result will be applied later to many concrete equations associated to the Schrödinger operator.

\begin{thm}
\label{thmestgen}
Let $\O\subset\R^N$ be a nonempty open subset of $\R^N,$ let $x_0\in\O,$ let $\rho_0>0,$ let $1\le p_1,\ldots,p_{n_1},q_1,\ldots,q_{n_2}<\infty,$ let
$F\in L^1_\loc(\O)$ be such that $F_{|\O\cap B(x_0,\rho_0)}\in L^2\big(\O\cap B(x_0,\rho_0)\big)$ and let
\begin{gather*}
f\in C\left(\bigcap\limits_{k=1}^{n_2}L^{q_k}_\loc(\O);\vsum_{j=1}^{n_1}L^{p_j^\p}_\loc(\O)\right).
\end{gather*}
Let $u\in H^1_\loc(\O)\cap L^{p_j}_\loc(\O)\cap L^{q_k}_\loc(\O),$ for any $(j,k)\in\li1,n_1\ri\times\li1,n_2\ri,$ be any solution to the complex-valued equation
\begin{gather}
\label{thmestgen1}
-\Delta u+f(u)=F, \text{ in } \Dr^\p(\O).
\end{gather}
If $\rho_0>\dist(x_0,\Gamma)$ then assume further that 
\begin{gather*}
f\in C\left(\bigcap\limits_{k=1}^{n_2}L^{q_k}(\O);\vsum_{j=1}^{n_1}L^{p_j^\p}(\O)\right), \; u\in H^1_0(\O), \\
u_{|\O\cap B(x_0,\rho_0)}\in L^{p_j}\big(\O\cap B(x_0,\rho_0)\big)\cap L^{q_k}\big(\O\cap B(x_0,\rho_0)\big),
\end{gather*}
for any $(j,k)\in\li1,n_1\ri\times\li1,n_2\ri.$ Set for every $\rho\in[0,\rho_0),$
\begin{gather}
\label{thmestgen2}
I(\rho)=\left|\dsp\int_{\O\cap\S(x_0,\rho)}u\ovl{\nabla u}.\frac{x-x_0}{|x-x_0|}\d\sigma\right|, \quad J(\rho)=\dsp\int_{\O\cap B(x_0,\rho)}|F(x)u(x)|\d x, \\
\label{thmestgen3}
w(\rho)=\int_{\O\cap\S(x_0,\rho)}u\ovl{\nabla u}.\frac{x-x_0}{|x-x_0|}\d\sigma, \quad
I_\Re(\rho)=\Re\big(w(\rho)\big), \quad I_\Im(\rho)=\Im\big(w(\rho)\big).
\end{gather}
Then we have,
\begin{gather}
\label{thmestgen4}
I,J,I_\Re,I_\Im\in C([0,\rho_0);\R), \\
\label{thmestgen5}
\|\nabla u\|_{L^2(\O\cap B(x_0,\rho))}^2+\Re\left(\:\vint_{\O\cap B(x_0,\rho)}f(u)\ovl u\d x\right)
=\Re\left(\:\vint_{\O\cap B(x_0,\rho)}F(x)\ovl{u(x)}\d x\right)+I_\Re(\rho), \\
\label{thmestgen6}
\Im\left(\:\vint_{\O\cap B(x_0,\rho)}f(u)\ovl u\d x\right)=\Im\left(\:\vint_{\O\cap B(x_0,\rho)}F(x)\ovl{u(x)}\d x\right)+I_\Im(\rho),
\end{gather}
for any $\rho\in[0,\rho_0).$
\end{thm}

\begin{rmk}
\label{rmkthmestgen}
One easily sees that if $\rho_0<\dist(x_0,\Gamma)$ then $I,J,I_\Re,I_\Im\in C([0,\rho_0];\R).$
\end{rmk}

\begin{exa}
\label{exathmestgen}
We give some functions $f$ for which Theorem~\ref{thmestgen} applies.
\begin{enumerate}
\item[]
\begin{enumerate}[$1)$]
\item
Typically, we apply Theorem~\ref{thmestgen} to
\begin{gather*}
f(u)=a|u|^{-(1-m)}u+bu+Vu,
\end{gather*}
with $(a,b)\in\C^2,$ $V\in L^\infty_\loc(\O)$ and $0<m<1.$ One easily checks that,
\begin{gather*}
f\in C\left(L^2_\loc(\O)\cap L^{m+1}_\loc(\O);L^2_\loc(\O)+L^\frac{m+1}{m}_\loc(\O)\right).
\end{gather*}
If in addition, $V\in L^\infty(\O)$ then one also has,
\begin{gather*}
f\in C\left(L^2(\O)\cap L^{m+1}(\O);L^2(\O)+L^\frac{m+1}{m}(\O)\right).
\end{gather*}
Let $z\in\C\setminus\{0\}.$ Since $\left||z|^{-(1-m)}z\right|=|z|^m,$ it is understood in the above example that $\left||z|^{-(1-m)}z\right|=0$ when $z=0.$
\item
\textbf{Hartree-Fock type equations.}
Let $V\in L^p(\R^N;\R)+L^\infty(\R^N;\R),$ with $\min\left\{1,\frac{N}{2}\right\}<p<\infty$ and let $W\in L^q(\R^N;\R)+L^\infty(\R^N;\R),$ with
$\min\left\{1,\frac{N}{4}\right\}<q<\infty.$ Set $r=\frac{2p}{p-1},$ $s=\frac{4q}{q-1},$
\begin{gather*}
E=L^2(\R^N)\cap L^4(\R^N)\cap L^r(\R^N)\cap L^s(\R^N),	\\
f(u)=Vu+(W\star|u|^2)u,
\end{gather*}
for any $u\in H^1(\R^N).$ Then $H^1(\R^N)\inj E$ with dense embedding and, by density of $\Dr(\R^N)$ in spaces $L^m(\R^N),$ for any
$m\in[1,\infty),$ we have
\begin{gather*}
E^\star=L^2(\R^N)+L^\frac43(\R^N)+L^{r^\p}(\R^N)+L^{s^\p}(\R^N),		\\
f\in C\big(E;E^\star\big),										\\
f\in C\big(H^1(\R^N);H^{-1}(\R^N)\big).
\end{gather*}
See Cazenave~\cite{MR2002047} (Proposition~1.1.3, Proposition~3.2.2, Remark~3.2.3, Proposition~3.2.9, Remark~3.2.10 and Example~3.2.11).
\end{enumerate}
\end{enumerate}
\end{exa}

\section{Proofs of the main results}
\label{proof}

Before proceeding to the proof of Theorems~\ref{thmsta} and \ref{thmstaF}, we recall the well-known Young's inequality. For any real $x\ge0,$ $y\ge0,$
$\la>1$ and $\eps>0,$ one has
\begin{gather}
\label{young}
xy\le\frac{1}{\la^\p}\eps^{\la^\p}x^{\la^\p}+\frac1\la\eps^{-\la}y^\la.
\end{gather}

\begin{vproof}{of Theorems~\ref{thmsta} and \ref{thmstaF}.}
We write $\rho_\star=\rho_0,$ for the proof of Theorem~\ref{thmsta} and $\rho_\star=\rho_1,$ for the proof of Theorem~\ref{thmstaF}. Let us introduce some notations. Let $\rho\in(0,\rho_\star).$ We set
\begin{gather*}
\begin{array}{lll}
E(\rho)=\|\nabla u \|_{L^2(B(x_0,\rho))}^2,		&	b(\rho)=\|u\|_{L^{m+1}(B(x_0,\rho))}^{m+1},	&	a(\rho)=\|u\|_{L^2(B(x_0,\rho))}^2,	\medskip\\
\theta=\frac{(1+m)+N(1-m)}{k}\in(0,1),		&	\ell=\frac{1}{\theta(1+m)},					&	\delta=\frac{k}{2(1+m)}.
\end{array}
\end{gather*}
We may assume that $u\in H^1\big(B(x_0,\rho_\star)\big).$ Indeed, the case $u\in H^1_\loc\big(B(x_0,\rho_\star)\big)$ can be treated by following the method in Bégout and D\'iaz~\cite{MR2876246} (see the end of Step~6, p.49, for Theorem~\ref{thmsta} and the end of Step~7, p.50, for
Theorem~\ref{thmstaF}. We now proceed with the proof in 3 steps. \\
{\bf Step~1.} $E\in W^{1,1}(0,\rho_\star),$ for a.e.~$\rho\in(0,\rho_\star),$ $E^\p(\rho)=\|\nabla u\|_{L^2(\S(x_0,\rho))}^2$ and
\begin{gather}
\label{proofthmsta1}
E(\rho)+b(\rho)\le\frac12
\left(K_1(\tau)\rho^{-(\nu-1)}E^\p(\rho)\right)^\frac{1}{2}\left(E(\rho)+b(\rho)\right)^\frac{\gamma(\tau)+1}{2}
+(L_1M)^2\|F\|_{L^2(B(x_0,\rho))}^2,
\end{gather}
where $K_1(\tau)=C(N,m)L_1^2M^2\max\left\{\rho_\star^{\nu-1},1\right\}\max\{b(\rho_\star)^{\mu(\tau)},b(\rho_\star)^{\eta(\tau)}\}$ and
$L_1=\max\left\{1,\frac{1}{L}\right\}.$ \\
By the first lines of Step~2, p.47, in Bégout and D\'iaz~\cite{MR2876246}, we only have to show \eqref{proofthmsta1}. Let $\rho\in(0,\rho_\star).$ We have to slightly modify the proof of Bégout and D\'iaz~\cite{MR2876246}. Indeed, since $F\in L^2,$ we need of the term $\|u\|_{L^2}^2.$ We have,
\begin{gather}
\label{proofthmsta2}
\left|\dsp\int_{\S(x_0,\rho)}u\ovl{\nabla u}.\frac{x-x_0}{|x-x_0|}\d\sigma\right|\le E^\p(\rho)^\frac{1}{2}\|u\|_{L^2(\S(x_0,\rho))}, \\
\label{proofthmsta3}
\|u\|_{L^2(\S(x_0,\rho))}\le C(N,m)\left(\|\nabla u\|_{L^2(B(x_0,\rho))}
+\rho^{-\delta}\|u\|_{L^{m+1}(B(x_0,\rho))}\right)^\theta\|u\|_{L^{m+1}(B(x_0,\rho))}^{1-\theta}.
\end{gather}
See Bégout and D\'iaz~\cite{MR2876246}: estimates (7.11), p.47, and (7.12), p.48. Putting together \eqref{thmsta1} (for Theorem~\ref{thmsta}), \eqref{thmstaF1} (for Theorem~\ref{thmstaF}), \eqref{proofthmsta2} and \eqref{proofthmsta3}, we obtain,
\begin{multline}
\label{proofthmsta4}
E(\rho)+b(\rho)+\kappa a(\rho)																					\\
\le CL_1ME^\p(\rho)^\frac{1}{2}\left(E(\rho)^\frac{1}{2}+\rho^{-\delta}b(\rho)^\frac{1}{m+1}\right)^\theta b(\rho)^\frac{1-\theta}{m+1}
+L_1M\vint_{B(x_0,\rho)}|F(x)u(x)|\d x,
\end{multline}
where $\kappa=0,$ in the case of Theorem~\ref{thmsta} and where $\kappa=1,$ in the case of Theorem~\ref{thmstaF}. In the case of
Theorem~\ref{thmstaF}, we apply~\eqref{young} with $x=|F|,$ $y=|u|,$ $\la=2$ and $\eps=\sqrt{L_1M},$ and we get
\begin{gather}
\label{proofthmsta5}
\vint_{B(x_0,\rho)}|F(x)u(x)|\d x\le\frac{L_1M}{2}\|F\|_{L^2(B(x_0,\rho))}^2+\frac{1}{2L_1M}a(\rho),
\end{gather}
for any $\rho\in(0,\rho_\star).$ Putting together \eqref{proofthmsta4} and \eqref{proofthmsta5}, we obtain for both theorems, for a.e. $\rho\in(0,\rho_\star),$
\begin{gather}
\label{proofthmsta6}
E(\rho)+b(\rho)\le
C_0L_1ME^\p(\rho)^\frac{1}{2}\left(E(\rho)^\frac{1}{2}+\rho^{-\delta}b(\rho)^\frac{1}{m+1}\right)^\theta b(\rho)^\frac{1-\theta}{m+1}
+(L_1M)^2\|F\|_{L^2(B(x_0,\rho))}^2.
\end{gather}
Let $\tau\in\left(\frac{m+1}{2},1\right]$ and let $\rho\in(0,\rho_\star).$ A straightforward calculation yields
\begin{align*}
     & \; \left(E(\rho)^\frac{1}{2}+\rho^{-\delta}b(\rho)^\frac{1}{m+1}\right)b(\rho)^\frac{1-\theta}{\theta(m+1)} \medskip \\
 =  & \; E(\rho)^\frac{1}{2}b(\rho)^\frac{1-\theta}{\theta(m+1)}+\rho^{-\delta}b(\rho)^\frac{1}{\theta(m+1)} \medskip \\
 =  & \; E(\rho)^\frac{1}{2}b(\rho)^{\tau(1-\theta)\ell}b(\rho)^{(1-\tau)(1-\theta)\ell}
            +\rho^{-\delta}b(\rho)^{\frac{1}{2}+\tau(1-\theta)\ell}b(\rho)^{\ell-\tau(1-\theta)\ell-\frac{1}{2}} \medskip \\
\le & \; 2\rho^{-\delta}\max\left\{\rho_\star^\delta,1\right\}K_2^2(\tau)^\frac{1}{2\theta}
            \left(E(\rho)+b(\rho)\right)^{\frac{1}{2}+\tau(1-\theta)\ell},
\end{align*}
where $K_2^2(\tau)=\max\{b(\rho_\star)^{\mu(\tau)},b(\rho_\star)^{\eta(\tau)}\},$ since $\frac{\mu(\tau)}{2\theta}=(1-\tau)(1-\theta)\ell$ and
$\frac{\eta(\tau)}{2\theta}=\ell-\tau(1-\theta)\ell-\frac{1}{2}.$ Hence~\eqref{proofthmsta1} follows from~\eqref{proofthmsta6} and the above estimate with $K_1(\tau)=16C_0^2L_1^2M^2K_2^2(\tau)\max\left\{\rho_\star^{\nu-1},1\right\},$ since $2\delta\theta=\nu-1$ and
$\theta\left(\frac{1}{2}+\tau(1-\theta)\ell\right)=\frac{\gamma(\tau)+1}{2}.$ \\
{\bf Step~2.} For any $\tau\in\left(\frac{m+1}{2},1\right]$ and for a.e.~$\rho\in(0,\rho_\star),$
\begin{gather*}
0\le E(\rho)^{1-\gamma(\tau)}\le K_1(\tau)\rho^{-(\nu-1)}E^\p(\rho)
+(2L_1M)^{2(1-\gamma(\tau))}\|F\|_{L^2(B(x_0,\rho))}^{2(1-\gamma(\tau))}.
\end{gather*}
Following Step~4, p.48, in Bégout and D\'iaz~\cite{MR2876246} but with Young's inequality~\eqref{young} applied with
$x=\frac12\left(K_1(\tau)\rho^{-(\nu-1)}E^\p(\rho)\right)^\frac{1}{2},$
$y=\left(E(\rho)+b(\rho)\right)^\frac{\gamma(\tau)+1}{2},$ $\la=\la(\tau)=\frac{2}{\gamma(\tau)+1}$ and
$\eps=\eps(\tau)=(\gamma(\tau)+1)^\frac{1}{\la(\tau)},$ Step~2 follows from the estimates
\begin{align*}
     &  \; E(\rho)+b(\rho) \\
\le & \; \frac12\left(K_1(\tau)\rho^{-(\nu-1)}E^\p(\rho)\right)^\frac{1}{2}\left(E(\rho)+b(\rho)\right)^\frac{\gamma(\tau)+1}{2}
															+(L_1M)^2\|F\|_{L^2(B(x_0,\rho))}^2,	\\
\le & \; \frac{C(\tau)}{2^\frac{\la(\tau)}{\la(\tau)-1}}\left(K_1(\tau)\rho^{-(\nu-1)}E^\p(\rho)\right)^\frac{1}{1-\gamma(\tau)}+\frac{1}{2}(E(\rho)+b(\rho))
															+(L_1M)^2\|F\|_{L^2(B(x_0,\rho))}^2,	\\
\le & \; \frac12\left(K_1(\tau)\rho^{-(\nu-1)}E^\p(\rho)\right)^\frac{1}{1-\gamma(\tau)}+\frac{1}{2}(E(\rho)+b(\rho))
															+(L_1M)^2\|F\|_{L^2(B(x_0,\rho))}^2,	\\
C(\tau)= & \; \frac{\la(\tau)-1}{\la(\tau)}\eps(\tau)^\frac{\la(\tau)}{\la(\tau)-1}
<\frac{\la\left(\frac{m+1}{2}\right)-1}{\la\left(\frac{m+1}{2}\right)}(\gamma(\tau)+1)^\frac{1}{\la(\tau)-1}
<\frac122^\frac{1}{\la(\tau)-1}<\frac122^\frac{\la(\tau)}{\la(\tau)-1}.
\end{align*}
{\bf Step~3.} Conclusion. \\
Now, following from Step~5 to Step~7, p.48--50, in Bégout and D\'iaz~\cite{MR2876246}, where estimate (7.16) therein has to be replaced with estimate of the above Step~2 and where  the mapping $\rho\longmapsto F(\rho)$ has to be replaced with the new function
$\rho\longmapsto(2L_1M)^{2(1-\gamma)}\|F\|_{L^2(B(x_0,\rho))}^{2(1-\gamma)},$ we prove Theorems~\ref{thmsta} and \ref{thmstaF}. This achieves the proof.
\medskip
\end{vproof}

\begin{vproof}{of Theorem~\ref{thmestgen}.}
If $\rho_0>\dist(x_0,\Gamma)$ then $u\in H^1_0(\O).$ So we may extend $u$ by $0$ on $\O^\c\cap B(x_0,\rho_0).$ Denoting $\wt u$ this extension, we have $\wt u\in H^1_0\big(\O\cup B(x_0,\rho_0)\big).$ We first consider the case where $\rho_0\neq\dist(x_0,\Gamma).$ We deal with
$\rho_0=\dist(x_0,\Gamma)$ at the end of the proof. It follows that $J\in C([0,\rho_0];\R)$ and by Cauchy-Schwarz's inequality, $I\in L^1(0,\rho_0).$ Thus, $I,J,I_\Re,I_\Im$ are defined almost everywhere on $(0,\rho_0).$ It follows from \eqref{thmestgen1} that,
\begin{gather}
\label{proofthmestgen1}
\langle\nabla u,\nabla\vphi\rangle_{\Dr^\p(\O),\Dr(\O)}+\langle f(u),\vphi\rangle_{\Dr^\p(\O),\Dr(\O)}=\langle F,\vphi\rangle_{\Dr^\p(\O),\Dr(\O)},
\end{gather}
for any $\vphi\in\Dr(\O).$ Let $\rho\in(0,\rho_0).$ For any $n\in\N,$ $n>\frac{1}{\rho},$ we define $\psi_n\in W^{1,\infty}(\R;\R)$ by
\begin{gather*}
\forall t\in\R,\; \psi_n(t)=
 \begin{cases}
  \: 1,			& \mbox{ if } |t|\in\left[0,\rho-\frac{1}{n}\right], \medskip \\
  \: n(\rho-|t|),	& \mbox{ if } |t|\in\left(\rho-\frac{1}{n},\rho\right), \medskip \\
  \; 0,			& \mbox{ if } |t|\in[\rho,\infty),
 \end{cases}
\end{gather*}
and we set $\wt\vphi_n(x)=\psi_n(|x-x_0|)\wt{u}(x)$ and $\vphi_n=\wt\vphi_{n|\O},$ for almost every $x\in\O\cup B(x_0,\rho_0).$ We easily check that for any $(j,k)\in\li1,n_1\ri\times\li1,n_2\ri,$
\begin{gather*}
\vphi_{n|{\O\cap B(x_0,\rho_0)}}
		\in H^1_0\big(\O\cap B(x_0,\rho_0)\big)\cap L^{p_j}\big(\O\cap B(x_0,\rho_0)\big)\cap L^{q_k}\big(\O\cap B(x_0,\rho_0)\big),	\\
\wt\vphi_n\in H^1_0\big(\O\cup B(x_0,\rho_0)\big)\cap L^{p_j}\big(\O\cup B(x_0,\rho_0)\big)\cap L^{q_k}\big(\O\cup B(x_0,\rho_0)\big),	\\
\vphi_n\in H^1_0(\O)\cap L^{p_j}(\O)\cap L^{q_k}(\O).
\end{gather*}
Then there exists $(\vphi_n^m)_{m\in\N}\subset\Dr(\O)$ such that for any $(n,m)\in\N^2,$ $\supp\vphi_n^m\subset\O\cap B(x_0,\rho_0)$ and
\begin{gather*}
\vphi_n^m\xrightarrow[m\tends\infty]{H^1_0(\O)\cap L^{p_j}(\O)\cap L^{q_k}(\O)}\vphi_n,
\end{gather*}
for any $(j,k)\in\li1,n_1\ri\times\li1,n_2\ri.$ Consequently, $\vphi=\vphi_n$ are admissible test functions in \eqref{proofthmestgen1}. We have,
\begin{align*}
	& \; \big\langle\nabla u,\nabla\vphi_n\big\rangle_{L^2(\O),L^2(\O)}
			=\langle\nabla\wt u,\nabla\wt\vphi_n\rangle_{L^2(\O\cup B(x_0,\rho_0)),L^2(\O\cup B(x_0,\rho_0))}				\\
   =	& \; \vint_{B(x_0,\rho)}\psi_n\big(|x-x_0|\big)|\nabla\wt u|^2\d x-n\,\Re\left(\;\vint_{\rho-\frac{1}{n}}^\rho\left(\:\vint_{\S(x_0,r)}
		\ovl{\wt {u}}\nabla\wt{u}.\frac{x-x_0}{|x-x_0|}\d\sigma\right)\d r\right),
\end{align*}
where we introduced the spherical coordinates $(r,\sigma)$ at the last line. We now let $n\nearrow\infty.$ Using the Lebesgue's dominated convergence Theorem and recalling that $I_\Re\in L^1((0,\rho_0);\R),$ we obtain
\begin{gather}
\label{proofthmestgen2}
\lim_{n\to\infty}\big\langle\nabla u,\nabla\vphi_n\big\rangle_{L^2(\O),L^2(\O)}=\|\nabla u\|^2_{L^2(\O\cap B(x_0,\rho))}-I_\Re(\rho).
\end{gather}
Proceeding as above but also with $\vphi=\vi\vphi_n,$ we get $\vlim_{n\to\infty}\big\langle\nabla u,\vi\nabla\vphi_n\big\rangle_{L^2,L^2}=-I_\Im(\rho)$ and
\begin{gather*}
	\begin{array}{ll}
		\vlim_{n\to\infty}\langle f(u),\vphi_n\rangle_{F^\star,E}=\Re\big(A(u)\big),
		&	\vlim_{n\to\infty}\langle f(u),\vi\vphi_n\rangle_{F^\star,E}=\Im\big(A(u)\big),	\medskip \\
		\vlim_{n\to\infty}\langle F,\vphi_n\rangle_{L^2,L^2}=\Re\big(B(u)\big),
		&	\vlim_{n\to\infty}\langle F,\vi\vphi_n\rangle_{L^2,L^2}=\Im\big(B(u)\big).
	\end{array}
\end{gather*}
where $E=\dsp\bigcap\limits_{j=1}^{n_2}L^{q_j}(\O),$ $F=\dsp\bigcap\limits_{j=1}^{n_1}L^{p_j}(\O),$ $A(u)=\dsp\int_{\O\cap B(x_0,\rho)}f(u)\,\ovl u\,\d x$ and
$B(u)=\dsp\int_{\O\cap B(x_0,\rho)}F(x)\ovl{u(x)}\d x.$ Estimates~\eqref{thmestgen5} and \eqref{thmestgen6} then follow from \eqref{proofthmestgen2} and these five last estimates. Since all terms in \eqref{thmestgen5} and \eqref{thmestgen6} are continuous on $[0,\rho_0],$ except eventually $I_\Re$ and $I_\Im,$ we deduce that $I_\Re$ and $I_\Im$ are continuous and \eqref{thmestgen5} and \eqref{thmestgen6} hold for any $\rho\in[0,\rho_0].$ The case $\rho_0=\dist(x_0,\Gamma)$ follows from the above proof applied with $\rho_0^n=\rho_0-\frac1n$ in place of $\rho_0$ and letting
$n\nearrow\infty.$
\medskip
\end{vproof}

\section{Application to the localization property to the case of Neumann boundary conditions}
\label{app}

In Bégout and D\'iaz~\cite{MR2876246}, the authors study the localization property for equation~\eqref{1} below with the homogeneous Dirichlet boundary condition (see, for instance, Theorem~3.5 in Bégout and D\'iaz~\cite{MR2876246}). In Theorem~\ref{thmcomneu1} below, we show that the same property holds with the homogeneous Neumann boundary condition. Before, we need to prove that solutions exist. This can be found in Bégout and D\'iaz~\cite{MR3315701}. Note that from Bégout and D\'iaz~\cite{MR2876246} to this paper, there was a slight change of notation. See Remark~\ref{com} below.

\begin{rmk}
\label{com}
In the context of the paper of Bégout and D\'iaz~\cite{MR2876246}, we can establish an existence result with the homogeneous Neumann boundary condition (instead of the homogeneous Dirichlet condition) and $F\in L^2(\O)$ \big(instead of $F\in L^\frac{m+1}{m}(\O)\big).$ In Bégout and
D\'iaz~\cite{MR2876246}, we introduced the set,
\begin{gather*}
\wt\A=\C\setminus\big\{z\in\C; \Re(z)=0 \text{ and } \Im(z)\le0\big\},
\end{gather*}
and assumed that $(\wt a,\wt b)\in\C^2$ satisfies,
\begin{gather}
\label{abtilde}
(\wt a,\wt b)\in\wt\A\times\wt\A \quad \text{ and } \quad
	\begin{cases}
		\Re(\wt a)\Re(\wt b)\ge0,													\medskip \\
		\text{ or }																\medskip \\
		\Re(\wt a)\Re(\wt b)<0 \; \text{ and } \; \Im(\wt b)>\dfrac{\Re(\wt b)}{\Re(\wt a)}\Im(\wt a),
	\end{cases}
\end{gather}
with possibly $\wt b=0,$ and we worked with
\begin{gather*}
-\vi\Delta u+\wt a|u|^{-(1-m)}u+\wt bu=\wt F.
\end{gather*}
But here in order to follow a closer notation with most of the works dealing with Schrödinger equations, we do not work any more with this equation but with,
\begin{gather*}
-\Delta u+a|u|^{-(1-m)}u+bu=F,
\end{gather*}
and $b\neq0.$ This means that we choose, $\wt a=\vi a,$ $\wt b=\vi b$ and $\wt F=\vi F.$ Then assumptions on $(a,b)$ are changed by the fact that for
$\wt z=\vi z,$
\begin{gather}
\label{sca1}
\Re(z)=\Re(-\vi\wt z)=\Im(\wt z),		\\
\label{sca2}
\Im(z)=\Im(-\vi\wt z)=-\Re(\wt z).
\end{gather}
It follows that the set $\wt\A$ and \eqref{abtilde} become,
\begin{gather}
\label{A}
\A=\C\setminus\big\{z\in\C; \Re(z)\le0 \text{ and } \Im(z)=0\big\},	
\end{gather}
\begin{gather}
\label{ab}
(a,b)\in\A\times\A \quad \text{ and } \quad
	\begin{cases}
		\Im(a)\Im(b)\ge0,										\medskip \\
		\text{ or }												\medskip \\
		\Im(a)\Im(b)<0 \; \text{ and } \; \Re(b)>\dfrac{\Im(b)}{\Im(a)}\Re(a).
	\end{cases}
\end{gather}
Obviously,
\begin{gather*}
\Big((\wt a,\wt b)\in\C^2 \text{ satisfies }\eqref{abtilde}\Big) \iff \Big((a,b)\in\C^2 \text{ satisfies }\eqref{ab}\Big).
\end{gather*}
\end{rmk}

\noindent
Assumptions~\eqref{ab} are made to prove the existence and the localization property of solutions to
\begin{gather}
\label{1}
-\Delta u+a|u|^{-(1-m)}u+bu=F, \text{ in } L^2(\O).
\end{gather}
For uniqueness, the hypotheses are the following (Theorem~2.10 in Bégout and D\'iaz~\cite{MR3315701}).

\begin{ass}[\textbf{Uniqueness}]
\label{ass}
Assume that $(a,b)\in\C^2$ satisfies one of the two following conditions.
\begin{enumerate}
 \item[]
  \begin{enumerate}[$1)$]
   \item
    $a\neq0,$ $\Re(a)\ge0$ and $\Re(a\ovl b)\ge0.$
   \item
    $b\neq0,$ $\Re(b)\ge0$ and $a=kb,$ for some $k\ge0.$
  \end{enumerate}
\end{enumerate}
\end{ass}

\noindent
A geometric interpretation of~\eqref{ab} and 1) of Assumption~\ref{ass} is given in Section~6 in Bégout and D\'iaz~\cite{MR2876246}, modulus a rotation in the complex plane. Now, we give some results about equation~\eqref{1} when $(a,b)\in\C^2$ satisfies \eqref{ab}.

\begin{cor}[\textbf{Neumann boundary conditions}]
\label{corexineu}
Let $\O$ be a nonempty bounded open subset of $\R^N$ having a $C^1$ boundary, let $\nu$ be the outward unit normal vector to $\Gamma,$ let $0<m<1$ and let $(a,b)\in\C^2$ satisfies~\eqref{ab}. For any $F\in L^2(\O),$ there exists at least one solution $u\in H^1(\O)$ to
\begin{gather}
	\label{corexineu1}
		\begin{cases}
			-\Delta u+a|u|^{-(1-m)}u+bu=F, \text{ in } L^2(\O),	\\
			\dfrac{\partial u}{\partial\nu}_{|\Gamma}=0.
		\end{cases}
\end{gather}
If furthermore $(a,b)$ satisfies Assumption~$\ref{ass}$ then the solution of~\eqref{corexineu1} is unique. Let $v\in H^1(\O)$ be any solution
to~\eqref{corexineu1}. Then $v\in  H^2_\loc(\O).$ In addition,
\begin{gather}
	\label{corexineu2}
		\|v\|_{H^1(\O)}\le M\|F\|_{L^2(\O)},
\end{gather}
where $M=M(|a|,|b|).$ Finally, if for some $\alpha\in(0,m],$ $F\in C^{0,\alpha}_\loc(\O)$ then $u\in C^{2,\alpha}_\loc(\O).$
\end{cor}

\begin{sym}
\label{sym}
If furthermore, for any $\vR\in SO_N(\R),$ $\vR\O=\O$ and if $F$ is spherically symmetric then we may construct a solution which is additionally spherically symmetric. For $N=1,$ this means that if $F$ is an even $($respectively, an odd$)$ function then $u$ is also an even $($respectively, an odd$)$ function.
\end{sym}

\noindent
Here and in what follows, $SO_N(\R)$ denotes the special orthogonal group of $\R^N.$ 

\begin{rmk}
\label{rmkexineu}
One easily checks that if $(a,b)\in\A^2$ satisfies $\Re(a)\ge0$ and $\Re(a\ovl b)\ge0$ then $(a,b)\in\C^2$ verifies~\eqref{ab}. In this case, uniqueness assumptions imply existence assumptions.
\end{rmk}

\begin{vproof}{of Corollary~\ref{corexineu} and Symmetry Property~\ref{sym}.}
The result comes from Bégout and D\'iaz~\cite{MR3315701}: Theorem~2.8 (existence and symmetry property), Theorem~2.10 (uniqueness), Theorem~2.9 \big(\textit{a priori} estimate~\eqref{corexineu2}\big) and Theorem~2.12 (local smoothness).
\medskip
\end{vproof}

\noindent
Concerning the support of solution of~\eqref{corexineu1} we have:

\begin{thm}
\label{thmcomneu1}
Let $\O$ be a nonempty bounded open subset of $\R^N$ having a $C^1$ boundary, let $0<m<1$ and let $(a,b)\in\C^2$ satisfies~\eqref{ab}. Then there exists $\eps_\star>0$ such that for any $0<\eps\le\eps_\star,$ there exists $\delta_0=\delta_0(\eps,|a|,|b|,N,m)>0$ satisfying the following property. Let $F\in L^2(\O)$ and let $u\in H^1(\O)$ be a solution to~\eqref{corexineu1}. If uniqueness holds for the problem~\eqref{corexineu1}\footnote{which is the case, for instance, if $(a,b)\in\C^2$ satisfies Assumption~\ref{ass}.}, $\supp F$ is a compact set and $\|F\|_{L^2(\O)}\le\delta_0$ then
$\supp u\subset K(\eps)\subset\O,$ where
\begin{gather*}
K(\eps)=\Big\{x\in\R^N;\; \exists y\in\supp F \; such \; that \; |x-y|\le\eps \Big\},
\end{gather*}
which is compact.
\end{thm}

\noindent
The proof relies on the following lemma.

\begin{lem}
\label{lemest}
Let $\O\subset\R^N$ be a nonempty open subset of $\R^N,$ let $0<m<1$ and let $(a,b)\in\C^2$ satisfies~\eqref{ab}. Let $F\in L^1_\loc(\O)$ and let
$u\in H^1_\loc(\O)$ be any solution to
\begin{gather}
	\label{lemest1}
		-\Delta u+a|u|^{-(1-m)}u+bu=F, \text{ in } \Dr^\p(\O).
\end{gather}
Then there exist two positive constants $L=L(|a|,|b|)$ and $M=M(|a|,|b|)$ satisfying the following property. Let $x_0\in\O$ and $\rho_\star>0.$ If
$F_{|\O\cap B(x_0,\rho_\star)}\in L^2\big(\O\cap B(x_0,\rho_\star)\big)$ then for any $\rho\in[0,\rho_\star),$
\begin{multline}
\label{lemest2}
\|\nabla u\|_{L^2(\O\cap B(x_0,\rho))}^2+L\|u\|_{L^{m+1}(\O\cap B(x_0,\rho))}^{m+1}+L\|u\|_{L^2(\O\cap B(x_0,\rho))}^2					\\
\le M\left(\left|\int_{\O\cap \S(x_0,\rho)}u\ovl{\nabla u}.\frac{x-x_0}{|x-x_0|}\d\sigma\right|+\int_{\O\cap B(x_0,\rho)}|F(x)u(x)|\d x\right),
\end{multline}
where it is additionally assumed that $u\in H^1_0(\O)$ if $\rho_\star>\dist(x_0,\Gamma).$
\end{lem}

\begin{proof*}
Let $x_0\in\O$ and let $\rho_\star>0.$ We set for every $\rho\in[0,\rho_\star),$
\begin{gather*}
I(\rho)=\left|\dsp\int_{\O\cap\S(x_0,\rho)}u\ovl{\nabla u}.\frac{x-x_0}{|x-x_0|}\d\sigma\right| \text{ and } J(\rho)=\dsp\int_{\O\cap B(x_0,\rho)}|F(x)u(x)|\d x.
\end{gather*}
It follows from Theorem~\ref{thmestgen} that $I,J\in C([0,\rho_\star);\R)$ and
\begin{gather}
\label{prooflemest1}
\left|\|\nabla u\|_{L^2(\O\cap B(x_0,\rho))}^2+\Re(a)\|u\|_{L^{m+1}(\O\cap B(x_0,\rho))}^{m+1}+\Re(b)\|u\|_{L^2(\O\cap B(x_0,\rho))}^2
\right|\le I(\rho)+J(\rho), \\
\label{prooflemest2}
\left|\Im(a)\|u\|_{L^{m+1}(\O\cap B(x_0,\rho))}^{m+1}+\Im(b)\|u\|_{L^2(\O\cap B(x_0,\rho))}^2\right|\le I(\rho)+J(\rho),
\end{gather}
for any $\rho\in[0,\rho_\star).$ Estimate \eqref{lemest2} then follows from \eqref{prooflemest1}, \eqref{prooflemest2} and Lemma~4.5 from Bégout and
D\'iaz~\cite{MR3315701} with $\delta=0.$
\medskip
\end{proof*}

\begin{vproof}{of Theorem~\ref{thmcomneu1}.}
Let $F\in L^2(\O)$ with $\supp F\subset\O$ and let $u\in H^1(\O)$ a solution to~\eqref{corexineu1} be given by Theorem~\ref{corexineu}. Set
$K=\supp F$ and
\begin{gather*}
\vO(\eps)=\Big\{x\in\R^N;\; \exists y\in K \text{ such that } |x-y|<\eps \Big\}.
\end{gather*}
Then $K(\eps)=\ovl{\vO(\eps)}.$ Let $\eps_\star>0$ be small enough to have $K(5\eps_\star)\subset\O$ and let $\eps\in(0,\eps_\star].$ Let $L$ and $M$ be given by Lemma~\ref{lemest} applied with $\rho_\star=2\eps.$ By Theorem~\ref{thmsta} and estimate~\eqref{corexineu2} in
Corollary~\ref{corexineu} above, there exists $\delta_0=\delta_0(\eps,|a|,|b|,N,m)>0$ such that if $\|F\|_{L^2(\O)}\le\delta_0$ then
$u_{\left|B(x_0,\eps)\right.}\equiv 0,$ for any $x_0\in\O$ such that $B(x_0,2\eps)\cap K=\emptyset$ and $B(x_0,2\eps)\subset\O.$ One easily sees that $B(x_0,2\eps)\cap K=\emptyset,$ for any $x_0\in\ovl{K(2\eps)^\c}\cap K(3\eps).$ We deduce that for any $x_0\in\ovl{K(2\eps)^\c}\cap K(3\eps),$ $u_{\left|B(x_0,\eps)\right.}\equiv 0.$ By compactness, there exist $n\in\N$ and $x_1,\ldots,x_n\in\ovl{K(2\eps)^\c}\cap K(3\eps)$ such that,
\begin{gather*}
\ovl{K(\eps)^\c}\cap\vO(4\eps)\subset\bigcup_{j=1}^nB(x_j,\eps)\subset\bigcup_{j=1}^nB(x_j,2\eps)\subset K(5\eps)\subset\O.
\end{gather*}
It follows that $u_{\left|K(\eps)^\c\cap\vO(4\eps)\right.}\equiv0.$ Let us define $\wt u$ in  $\O$ by,
\begin{gather*}
\wt u=
	\begin{cases}
		u,	&	\text{in } \vO(2\eps),			\\
		0,	&	\text{in } \O\setminus\vO(2\eps).
	\end{cases}
\end{gather*}
It follows that $\supp\wt u\subset K(\eps)$ and $\wt u\in H^1_0(\O)$ is a solution to~\eqref{corexineu1}. By uniqueness assumption, $\wt u=u$ so that
$\supp u\subset K(\eps)\subset\O,$ which is the desired result.
\medskip
\end{vproof}

\begin{rmk}
\label{rmkfinal}
In Bégout and D\'iaz~\cite{MR2876246}, the authors study existence, uniqueness, smoothness and localization property for the equations~\eqref{1} with an external source $F$ belonging to $L^\frac{m+1}{m}(\O)$ with $0<m<1$ (see, for instance, Theorem~3.5 in Bégout and D\'iaz~\cite{MR2876246}). Below, we explain how the same results hold true with the weaker assumption $F\in L^2(\O).$ Indeed, when $|\O|<\infty$ and $0<m<1,$ $L^\frac{m+1}{m}(\O)\inj L^2(\O)$ and $L^\frac{m+1}{m}(\O)\neq L^2(\O).$ Results of existence can be found in Bégout and D\'iaz~\cite{MR3315701} jointly to some others additional results. Hypotheses on $(a,b)\in\C^2$ are the same as in Bégout and D\'iaz~\cite{MR2876246}, except we have to require $b\neq0.$ Note that from Bégout and D\'iaz~\cite{MR2876246} to the present paper, there was a change of notation. See Remark~\ref{com} for precision. Throughout this remark, equation \eqref{1} with homogeneous Dirichlet boundary condition are considered and $F$ is always assumed to belong in $L^2(\O)$ (instead of $L^\frac{m+1}{m}(\O)$ in Bégout and D\'iaz~\cite{MR2876246}) and assumptions on $(a,b)$ are \eqref{ab} and Assumption~\ref{ass}, instead of (2.2) and (2.3) in Bégout and D\'iaz~\cite{MR2876246}.
\medskip
\\
Analogous results to Theorems~4.1, 4.4 and Corollary~5.3 of Bégout and D\'iaz~\cite{MR2876246} can be easily adapted. Indeed, by Bégout and
D\'iaz~\cite{MR3315701} (Theorems~2.8, 2.9, 2.10 and 2.12), these results hold but with $u\in H^2_\loc(\O)$ and
\begin{gather}
\label{corexidir2}
\|u\|_{H^1(\O)}^2+\|u\|_{L^{m+1}(\O)}^{m+1}\le M\|F\|_{L^2(\O)}^2,
\end{gather}
instead of $u\in W^{2,\frac{m+1}{m}}_\loc(\O),$ (4.1) and (4.2) in Bégout and D\'iaz~\cite{MR2876246}. Concerning the localization property, Theorems~3.1 and 3.5 in Bégout and D\'iaz~\cite{MR2876246} still hold true but with $F\in L^2(\O)$ and
\begin{gather}
\label{corF1}
\forall\rho\in(0,\rho_1), \;
\|F\|_{L^2(\O\cap B(x_0,\rho))}^2\le\eps_\star(\rho-\rho_0)_+^p,
\end{gather}
instead of (3.1) in Bégout and D\'iaz~\cite{MR2876246}. The proofs are essentially the same where we use Lemma~\ref{lemest} and \eqref{corexidir2} above instead of (4.1) in Bégout and D\'iaz~\cite{MR2876246}. It follows that Theorems~1.1 and 1.2 in Bégout and D\'iaz~\cite{MR2876246} can be easily adapted with the obvious modifications.
\end{rmk}

\end{document}